\documentstyle[11pt,amscd,amssymb]{amsart}
\newcommand{\N}{{\Bbb N}}
\newcommand{\Ztwo}{{\Bbb Z}^2}
\newcommand{\Ztwol}{{\Bbb Z}^2_l}
\newcommand{\e}{\mbox{e}}
\newcommand{\E}{{\Bbb E}}
\renewcommand{\P}{{\Bbb P}}
\newcommand{\R}{{\Bbb R}}
\newcommand{\df}{\ \stackrel{\Delta}{=}\ }

\newcommand{\dr}{\text{d}}
\begingroup 
\newtheorem{thm}{Theorem}
\newtheorem{lem}[thm]{Lemma}         
\newtheorem{prop}[thm]{Proposition}  
\endgroup

\theoremstyle{definition}

\theoremstyle{remark}
\newtheorem{rem}[thm]{Remark}        

\begin{document}
\title[\,]{A note on the decay of correlations 
 Under $\delta-$Pinning}
\author{Dmitry Ioffe}
\address{
WIAS, Mohrenstr.~39   
D-10117 Berlin, Germany  and  Faculty of Industrial Engineering,
Technion, Haifa 32000, Israel}
\author{ Yvan Velenik}
\address{Fachbereich Mathematik, Sekt.~MA~7-4, TU-Berlin, 
Strasse~des~17 Juni 136, D-10623 Berlin, Germany}
\date{\today}
\begin{abstract}
We prove that for a class of massless $\nabla\phi$ interface models on $\Ztwo$
an introduction of an arbitrary small  pinning self-potential  leads to
exponential decay of correlation, or, in other words, to  creation of mass.
\end{abstract}
\maketitle

\vskip 0.1in 

In this note we study a family of effective interface models over  $\Ztwo$ with
the formal Hamiltonian ${\cal H}$ given by
\begin{equation}
\label{hamiltonian}
{\cal H}(\phi )~=~\sum_{i\sim j} V(\phi_i -\phi_j ),
\end{equation}
where the summation is over all nearest neighbours $i\sim j$ of $\Ztwo$, and
the following two assumptions are made on the interaction potential $V$:
\begin{itemize}
\item $V$ is even and smooth
\item There exists a constant $c_V \geq 1$, such that
\begin{equation}
\label{V}
\frac1{c_V}~\leq~V^{\prime\prime} (t)~\leq~c_V\qquad\qquad\forall
t\in\R .
\end{equation}
\end{itemize}
\begin{rem}
No further assumptions on $c_V$ are made, and,  in fact, we expect that the
results of the paper remain true if only the lower bound in \eqref{V} is
assumed. Also,  though we do not stipulate it explicitly at each particular
instance, the values of {\bf all} the positive constants we use below depend on
$c_V$. 
\end{rem}

Given a set $A\subset\Ztwo$ with a finite complement $A^c\df\Ztwo\setminus A$,
we use $\P_A$ to denote the finite volume Gibbs measure on  $\Omega_A\df
\R^{A^c}$ with the   Hamiltonian ${\cal H}$ and zero boundary conditions on
$A$;
\begin{equation}
\label{PA}
\P_A (\dr \phi )~=~\frac1{{\bf Z}(A)}\e^{-{\cal H}(\phi )}\prod_{i\in
  A^c}
\dr h_i\prod_{j\in A}\delta_0 (\dr h_j ) .
\end{equation}
It is well known that $\P_A$ delocalizes as $A^c\nearrow\Ztwo$; maybe the
easiest way to see this is to use the reverse Brascamp-Lieb inequality
\cite{DGI} which implies that the variance of $\phi_0$ under $\P_A$ dominates the
corresponding Gaussian variance. If, however,
an, essentially arbitrary small, pinning self-potential is added to ${\cal H}$,
then the situations radically changes, and the infinite volume Gibbs state
exists in the usual sense. This phenomenon has  been first worked out in the
Gaussian case ($c_V =1$) in \cite{DMRR}. Our main reference \cite{DV} contains
a proof of the localization for a fairly general class of interactions and
self-potentials. In this note we prove that in the case of the family of random
interfaces as in \eqref{hamiltonian}, the delocalization/localization
transition is sharp in the sense that it always comes together with the
exponential decay of correlations, or, using the language of a more physically
oriented literature, with the creation of mass.

For simplicity, but also in order to give a cleaner exposition of  otherwise
more general renormalization ideas behind the proof, we consider here only the
case of the so called $\delta$-pinning, thereby generalizing recent results of
\cite{BB} on purely Gaussian fields (that is again $c_V=1$):

Given a box $\Lambda_N\df [-N,...,N]^2\subset\Ztwo$ and a number  $J\in\R$
(which characterizes the strength of the pinning) we define the following
measure $\hat{\P}_N$ on $\R^{\Lambda_N}$:

\begin{equation}
\label{PNhat}
\hat{\P}_N (\dr \phi )~=~\frac1{\hat{{\bf Z}}_N}\e^{-{\cal H}(\phi )}
\prod_{i\in\Lambda_N}\left(\dr \phi_i +\e^J\delta_0(\dr \phi_i
)\right)
\prod_{j\in\Ztwo\setminus\Lambda_N}\delta_0 (\dr \phi_j ) .
\end{equation}

Notice that the case $J= -\infty$ corresponds to the original measure on
$\R^{\Lambda_N}$ with the Hamiltonian \eqref{hamiltonian}, which delocalizes as
$N\to\infty$.
\begin{lem}
\label{Lemma1}
For every $J\in\R$ there exists an exponent (mass) $m=m(J) >0$ and a 
constant $c_1 = c_1 (J)<\infty$, such that
\begin{equation}
\label{mass}
{\Bbb C}\text{\rm ov}_{\hat{\P}_N}\big( \phi_i
;\phi_j\big)~\leq~c_1\e^{-m\|i-j\|}
\end{equation}
uniformly in $N$ and in $i,j\in\Ztwo$.
\end{lem}

Of course, there is nothing to prove  if either $i$ or $j$ lies outside of
$\Lambda_N$. In fact, the sub-index $N$ is superfluous - all the estimates we
use and obtain simply do not depend on a particular $\Lambda_N$, and the only
reason we need it is to make the definitions mathematically meaningful. From
now on we shall drop the sub-index $N$ from the notation.
\vskip 0.1in
A right way to think about \eqref{PNhat} is as of the joint distribution of the
field of random interface heights $\{\phi_i\}_{i\in\Ztwo}$  and the random
``dry'' set ${\cal A}$;
$$
{\cal A}~\df~\left\{ i\in\Ztwo :
 \phi_i =0\right\} .
$$
Integrating out all the height variables $\phi$ in \eqref{PNhat} we arrive to
the following probability distribution for ${\cal A}$;
\begin{equation}
\label{dryset}
\hat{\P}\left( {\cal A}=A\right)~\df~\rho (A)~=~\frac1{\hat{\bf Z}}
\e^{J|A|}{\bf Z} (A)~=~\frac{\e^{J|A|}{\bf Z} (A)}
{\sum_{D}\e^{J|D|}{\bf Z} (D)} ,
\end{equation}
where the partition function ${\bf Z}(A)$ is the same as in \eqref{PA}.

Using the probabilistic weights $\{\rho (A)\}$ one can rewrite $\hat{\P}$ as
 the convex combination,
\begin{equation}
\label{Phatconv}
\hat{\P}(\cdot )~=~\sum_A\rho (A)\P_A (\cdot ) .
\end{equation}
Since under each $\P_A$ the distribution of $\phi_i$ is symmetric for every
$i\in\Ztwo$, this gives rise to the following decomposition  of the
covariances:
\begin{equation}
\label{Covconv}
{\Bbb C}\text{ov}_{\hat{\P}}\big( \phi_i ;\phi_j\big)~=~
\sum_{A}\rho (A)\langle \phi_i ;\phi_j\rangle_{A} .
\end{equation}
 
At this point we shall utilize the random walk representation of $\langle
\phi_i ;\phi_j\rangle_{A}$ which has been first developed in the PDE context in
\cite{HS}. We follow the approach of \cite{DGI}, where the
Helffer-Sj\"{o}strand representation was put on the probabilistic tracks:

One constructs a stochastic process $\big(\Phi (t) ,X(t)\big)$,  where:
\begin{itemize}
\item $\Phi (\cdot ) $ is a diffusion on $\R^{A^c}$ with the  invariant measure
$\P_A$.
\item Given a trajectory $\phi (\cdot )$ of the process $\Phi$,  $X(t)$ is an,
in general  inhomogeneous, transient random walk on $A^c\cup\partial
A^c\subset\Ztwo$ with the life-time
$$
\tau_A~\df~\inf\{t:~X(t)\in A\},
$$
and the time-dependent jump rates
\begin{equation}
\label{rates}
a(i,j;t)~=~\left\{
\begin{split}
&V^{\prime\prime}\big(\phi_i (t) -\phi_j (t)\big),\qquad \text{if}\
i\sim j\\
&0,\qquad\qquad\qquad\ \qquad \text{otherwise}
\end{split}
\right.
\end{equation}
\end{itemize}

Let us use ${\cal E}_{i,\phi}^{A}$ to denote the law of  $\big( X(t) , \Phi
(t) \big)$ starting from the point  $(i,\phi )\in A^c\times \R^{A^c}$. Then
(\cite{HS},\cite{DGI}),
\begin{equation}
\label{RWrepr}
\langle\phi_i ,\phi_j\rangle_A~=~\left\langle 
{\cal E}_{i,\phi}^{A}\int_{0}^{\tau_A}{\Bbb I}_{\{X(s)=j\}}\dr
s\right\rangle_A.
\end{equation}

Substituting the latter expression into \eqref{Covconv},
\begin{equation}
\label{basic}
{\Bbb C}\text{ov}_{\hat{\P}}\big( \phi_i ;\phi_j\big)~=~
\sum_A\rho (A)
\left\langle 
{\cal E}_{i,\phi}^{A}\int_{0}^{\tau_A}{\Bbb I}_{\{X(s)=j\}}\dr
s\right\rangle_A.
\end{equation}

It is very easy now to explain the logic behind the proof of 
Lemma~\ref{Lemma1}: The expression
$$
{\cal E}_{i,\phi}^{A}\int_{0}^{\tau_A}{\Bbb I}_{\{X(s)=j\}}\dr s
$$
describes the time spent by the random walk $X(\cdot )$ starting at  $i$ in the
site $j$ before being killed upon entering the dry set $A$ which,  for the
purpose, could be considered as a random killing obstacle. In order to prove
that this time is exponentially (in $\| i-j\|$) small one  needs an appropriate
density estimate on $A$ and a certain path-wise control on the exit
distributions of $X (\cdot )$. In the Gaussian case considered in \cite{BB},
$X(\cdot )$ happens to be just the simple  random walk on $\Ztwo$ which is
completely decoupled from the  diffusion part $\Phi (\cdot )$, and, thus,
behaving independently of  $A$ and the initial condition $\phi\in\R^{A^c}$.
This lead in \cite{BB} to a resummation argument, which substantially
facilitated the matter. One of the main difficulties in the non-Gaussian case 
we consider here is the dependence of the distribution of $X(\cdot  )$ on the
realization of the dry set $A$ and on the sample path of the  diffusion $\Phi$.
We still have very little to say about this dependence. However, due to the
basic assumption \eqref{V} on the  interaction potential $V$, the jump rates
$a(i,j;t)$ in \eqref{rates} are uniformly bounded above and below:
\begin{equation}
\label{jumprates}
\frac1{c_V}~\leq a(i,j ;t)~\leq~c_V .
\end{equation}
In particular one always has a rough control over probabilities of hitting
distributions. For example, if the random walk $X$ enters a box ${\bf B}_l$ of
linear size $l$ which is known to contain a dry site; it would be convenient to
call such a box ``dirty'', then the  probability that $X$ hits this site (and
consequently dies there) before leaving ${\bf B}_l$ should be bounded below by
some positive number $p=p(l)>0$. Thus if the realisation $A$ of the random dry
set ${\cal A}$ is such, that on its way from $i$ to $j$ the walk $X$ cannot
avoid visiting less than $\epsilon \| i-j\|$ disjoint dirty $l$-boxes, the
probability that it eventually reaches $j$ before being killed should be
bounded above by something like
$$
\left( 1- p(l) \right)^{\epsilon \| i-j\|} .
$$
Proposition~\ref{Prop2} below makes this computation precise.

The crux of the matter, however, is to ensure that on a certain  finite
$l$-scale the density of the dirty $l$-boxes is so high, that only with
exponentially small probabilities the realization $A$ of ${\cal A}$ enables an
$\epsilon$-clean passage from $i$ to $j$. A statement of this  sort  is given
in Proposition~\ref{Prop1}. 
\vskip 0.1in
\noindent
Once the renormalization approach sketched above is accepted as the strategy of
the proof, the first drive of an associative thinking is to try to compare the
distribution of ${\cal A}$ on different $l$-scales with, say, independent
Bernoulli percolation or other known models with controllable decay of
connectivities. This we have tried and failed, and, at least in the case of
$\Ztwo$, such a comparison is unlikely.

The relevant statistical properties of the random dry set ${\cal A}$  on
various finite length scales are captured in the following estimate which
generalizes the key Proposition~5.1 in \cite{DV}
\begin{thm}
\label{Main}
For each $J\in \R$ there exists a number $R =R (J)<\infty$ and
 exponent $\nu =\nu (J) >0$, such that whenever a
finite set $B\subset\Ztwo$ admits a decomposition
\begin{equation}
\label{Bdecomp}
B~=~\bigvee_{l=1}^n B_l
\end{equation}
into connected disjoint components $B_1,...,B_n$ with
\begin{equation}
\label{diam}
\text{diam}\big( B_l\big)~\geq~R;\qquad\qquad l=1,...,n ,
\end{equation}
the following exponential upper bound on having all of $B$ ``clean of 
dry points'' 
holds:
\begin{equation}
\label{mainbound}
\sum_{A\cap B =\emptyset}\rho (A)~\leq~\e^{-\nu |B|}.
\end{equation}
\end{thm}

We relegate the proof of Theorem~\ref{Main} to the end of the 
paper, and, assuming for the moment its validity, directly proceed to
the proof of the mass-generation claim of Lemma~\ref{Lemma1}.

\vskip 0.2in

\noindent
{\it Proof of Lemma~\ref{Lemma1}:\ }\  The number $R=R(J)$ which  appears in
the basic Theorem~\ref{Main} sets up the stage for the finite scale 
renormalization analysis of the random dry set ${\cal A}$.  Let us pick a
number $l>R;\ l\in\N ,$ and define the renormalized  lattice
$$
\Ztwol~\df~(2l+1)\Ztwo .
$$
To distinguish between the sets on the original lattice $\Ztwo$ and  those on
the renormalized one $\Ztwol$ we shall always mark the latter by the
super-index $l$. For example ${\bf B}^l (x,r)$  stands for the $\Ztwol$ lattice
box centered at $x\in\Ztwol$;
$$
{\bf B}^l (x,r)~\df~\left\{ y\in\Ztwol :\ \ \|x-y\|\leq lr\right\} .
$$
Let us define $\Gamma^l (r)$ as the set of all $\Ztwol$-nearest neighbour
lattice paths leading from the origin to the boundary $\partial {\bf B}^l
(x,r)$. With each $\gamma^l\in \Gamma^l (r)$ we associate a connected chain
$\tilde{\gamma}^l$ of $l$-blocks on the original lattice $\Ztwo$;
$$
\tilde{\gamma}^l~\df~\bigcup_{x\in\gamma^l} {\bf B} (x,l) .
$$
Let us fix a  number $\epsilon \in (0,1)$. We say that a path
$\gamma^l\in\Gamma^l$ is $(r,\epsilon )$-clean in $A\subset\Ztwo$, if
$$
\#\left\{x\in\gamma^l :\ \ {\bf B} (x,l)\cap A \neq\emptyset\right\}~<~
\epsilon r .
$$
Similarly, we say that a set $A\subset\Ztwo$ is $(r,\epsilon )$-clean if there
exists a path $\gamma^l\in\Gamma^l (r)$ which is  $(r,\epsilon )$-clean for
$A$. Otherwise, we shall call $A$  $(r,\epsilon)-$dirty.

\begin{prop}
\label{Prop1} For each $\epsilon\in (0,1)$  there exist a number $l_0 =l_0
(\epsilon, J)<\infty$ and a radius $r_0 =r_0 (\epsilon )$,  such that for every
choice of $l\geq l_0$;
$$
\sum_{A~\text{is}~(r,\epsilon )-\text{clean}}\rho\left( 
A\right)~\leq~
\e^{-c_2 (\epsilon ,l)r} ,
$$
uniformly in $r\geq r_0$, where $c_2 (\epsilon ,l )$  diverges (as $l^2$) with
$l$.
\end{prop}
\vskip 0.1in
\noindent
{\it Proof:}\ \ 
The condition on $r_0(\epsilon)$ is a semantic one - the only thing  we want is
to ensure that $r >[\epsilon r]$.

 Let us estimate the probability of  the event  $\{ A~\text{is}~(r,\epsilon
)-\text{clean}\}$ as follows:
\begin{equation}
\label{epsilonclean}
\sum_{A\,\text{is}\,(r,\epsilon )-\text{clean}}\rho\left( 
A\right)~\leq~
\sum_{k=r}^{\infty}\sum_{\gamma^l\in\Gamma^l :|\gamma^l |=k}\ \ 
\sum_{A\,:\gamma_l\,\text{is}\,(r,\epsilon)-\text{clean in}\,A}\rho 
(A) .
\end{equation}
Each path $\gamma^l =(0,x_1,...,x_k);\ \gamma^l\in\Gamma^l$, which is
$(r,\epsilon )$-clean in $A$ contains at most $[\epsilon r]$  vertices
$x_{i_1},...,x_{i_M};\ M\leq [\epsilon r]$, such that the  corresponding
$l-$blocks have a non-empty intersection with $A$;
$$
{\bf B} (x_i,l)\cap A~\neq~\emptyset ;\qquad i=1,...,M .
$$
Whatever happens, for a path $\gamma^l$ of length $k$ there are at most  $2^k$
(in fact much less due to the restriction $M\leq [\epsilon r]$) possible ways
to choose  a sub-family $\tilde{\gamma}_{\text{dirty}}^l$;
$$
\tilde{\gamma}_{\text{dirty}}^l~\df~\bigcup_{i=1}^{M}{\bf B}(x_i ,l),
$$
 of ``dirty'' block along $\tilde{\gamma}^l$. On the other hand, fixing both
$\tilde{\gamma}^l$ and its  ``dirty part'' $\tilde{\gamma}^l_{\text{dirty}}$, 
we can use Theorem~\ref{Main} to obtain
\begin{equation}
\label{cleanfactor}
\sum_{A\cap\tilde{\gamma}^l\setminus\tilde{\gamma}^l_{\text{dirty}}=
\emptyset}\rho (A)~\leq~
\exp\{-\nu |
\tilde{\gamma}^l\setminus\tilde{\gamma}^l_{\text{dirty}}|\}~\leq~
\e^{-\nu (k-[\epsilon r])l^2}
\end{equation}
We, thus, conclude, 
that for any $k\geq r$ and for each 
$\gamma^l\in\Gamma^l$ with $|\gamma^l | =k$,
\begin{equation*}
\sum_{A:\,\gamma_l~\text{is}\,(r,\epsilon)-\text{clean in}\,A}\rho
(A)~\leq~
\e^{-\nu (J) l^2 (k-[\epsilon r]) +k\log 2} .
\end{equation*}
Using the above estimate together with  the trivial bound;
$$
\#\left\{ \gamma^l\in\Gamma^l :\ |\gamma^l |=k\right\}~\leq~4^k,
$$
 to perform the summation in \eqref{epsilonclean} we arrive at the claim of
Proposition~\ref{Prop1}.
\qed
\vskip 0.1in

\noindent
Nothing in the above argument depends on the fact that the box  ${\bf B}(0,rl)$
is centered at the origin. Without any loss of generality we shall prove
\eqref{mass} only for the case $i=0$. 

Let us fix $l$ and $\epsilon$ as in the statement of  Proposition~\ref{Prop1}.
For each $j$ with $\|j\|>rl$ we use  \eqref{basic} and estimate:
\begin{equation}
\label{covdecomp}
\begin{split}
{\Bbb C}\text{ov}\big( \phi_0 ;\phi_j\big)~\leq~
&\sum_{A~\text{is}~(r,\epsilon )-\text{clean}}\rho\left( A\right)\\
&+~\sum_{A~\text{is}~(r,\epsilon )-\text{dirty}}
\rho\left( A\right)\max_{\phi}{\cal E}_{0,\phi}^A\int_0^{\tau_A}
{\Bbb I}_{\{X(s)=j\}}\dr s .
\end{split}
\end{equation}

The first term in \eqref{covdecomp} has been just estimated in 
Proposition~\ref{Prop1}. Let us use $\tau_{rl}$ to denote the exit  time from
${\bf B} (0,rl)$. The second term in \eqref{covdecomp} could be further 
bounded above as
\begin{equation}
\label{second}
\max_{A~\text{is}~(r,\epsilon)-\text{dirty}}\max_{\phi}
 {\cal E}^A_{0,\phi} \big(\tau_A >\tau_{rl}\big)
\sum_B\rho(B)\max_{\psi}{\cal E}^B_{j,\psi}\int_0^{\tau_B} 
{\Bbb I}_{\{X(s)=j\}}\dr s.
\end{equation}

It is convenient to estimate the above expression in a complete generality of
time dependent random walks with bounded jump rates $a(i,j;t)$:

Let $X (t)$ be the time-inhomogeneous Markov process with the transition rates
as in \eqref{jumprates}. It is always possible to  homogenize it, and to
consider
$$
\tilde{X} (t)~\df~(X(t) ,t ) .
$$
We shall use $\tilde{\E}_{(i,t)}$ to denote the law of $\tilde{X}$ with the
space-time starting point $(i,t)\in\Ztwo\times\R$.

The ${\bf B}(0,rl)$ box is decomposed to the disjoint union of  sub-blocks on
the $l$-scale as:
$$
{\bf B}(0,rl)~=~\cup_{x\in {\bf B}^l (0,r)}{\bf B}(x ,l) .
$$
To a generic point $i\in\ {\bf B} (0,rl)$ we associate an $l$-block ${\bf B}_l
(i)$ according to the following rule:
$$
{\bf B}_l (i)~=~{\bf B}(x ,l)\ \ \text{if}\ i\in {\bf B}(x ,l)\ 
\text{for some}\ x\in\Ztwol .
$$
Given a $(r,\epsilon )$-dirty set $A\subset\Ztwo$, let us call a  block ${\bf
B}(x ,l);\ x\in\Ztwol$, dirty if
$$
{\bf B}(x ,l)\cap A\neq\emptyset .
$$

We introduce now the following  family of stopping times for the process
$\tilde{X} (t)$:
\begin{equation*}
\begin{split}
&T_1 =\inf_{t\geq 0}\{ {\bf B}_l (X(t ))\ \text{is dirty}\}. \\
&\,\\
&S_1 =\inf_{t\geq T_1}\{ {\bf B}_l (X(t))\neq {\bf B}_l
  (X(T_1))\} . \\
&\,\\
&T_2 =\inf_{t\geq S_1}\{{\bf B}_l (X(t))\ \ \text{is dirty}\}\\
&............................\\
&S_n =\inf_{t\geq T_n}\{ {\bf B}_l (X(t))\neq {\bf B}_l
  (X(T_n))\}. \\
&........................... \\
\end{split}
\end{equation*}

The condition of $A$ being $(r,\epsilon )$-dirty is readily translatable under
$\P_A$  to the sure event
$$
\left\{\tau_{rl}~>~T_{\epsilon r}\right\} .
$$
Consequently, if, as before, we use $\tau_A$ to denote the hitting time of the
set $A$ ,
$$
\tilde{\P}_{(0,0)} (\tau_A >\tau_{rl})~\leq~\tilde{\P}_{(0,0)}
 (\tau_A > T_{\epsilon r} )~=~
\tilde{\E}_{(0,0)}\tilde{\E}_{\tilde{X} (T_1 )} {\Bbb I}_{\tau_A >
    S_1}...
\tilde{\E}_{\tilde{X}(T_{\epsilon r})}{\Bbb I}_{\tau_A >
    S_{\epsilon r}} .
$$
We claim that each of the $\epsilon r$ terms in the above product admits  an
upper bound of the form
\begin{equation}
\label{dirtybound}
1~-~\left(\frac1{3c_V^2 +1}\right)^{2l} .
\end{equation}
uniformly in all Markov chains with bounded rates condition \eqref{jumprates} 
and (which is the same)  in all possible values of above stopping times.

Indeed let ${\bf B}_l$ be a box of side length $l$, and $i,k\in {\bf B}_l$ .
Then one strategy for a random walk starting at $i$ to hit $k$ before leaving
${\bf B}_l$ is to march to $k$ directly  along some prescribed unambiguous
trajectory, say first horizontally and then vertically. Clearly if one pulls
down the  rates along such a trajectory to the minimum value $1/c_V$ and pushes
the rates leading out of this trajectory to the maximal value $c_V$, then the
probability to follow the trajectory itself only decreases, but to an exactly
computable value
$$
\left(\frac1{3c_V^2 +1}\right)^{\|i-k\|} ,
$$
where the power $\|i-k\|$, of course, corresponds to the number of steps along
the trajectory. Hence \eqref{dirtybound}.

As a result:
\begin{prop}
\label{Prop2}
Uniformly in $r$ and in $(r,\epsilon)$-dirty sets $A$,
$$
\max_{\phi} {\cal E}_{0,\phi}^A \left(\tau_A >\tau_{rl}\right)
~\leq~\e^{-c_3 rl} .
$$
\end{prop}

\qed
\vskip 0.1in
\noindent
Finally, 
\begin{equation}
\label{selfterm}
\begin{split}
\sum_B\rho(B)\max_{\phi}&{\cal E}^B_{j,\phi}\int_0^{\tau_B} 
{\Bbb I}_{\{X(s)=j\}}\dr s\\
&=~\sum_{k=1}^{\infty}\sum_{B:\text{d}(j,B)=k}
\rho(B)\max_{\phi}{\cal E}^B_{j,\phi}\int_0^{\tau_B} 
{\Bbb I}_{\{X(s)=j\}}\dr s ,
\end{split}
\end{equation}
where $\text{d}(j,B)\df\inf\{ \| j-i\| :\ i\in B\}$.

Proceeding as in the proof of Proposition~\ref{Prop2}, we readily obtain that
there exists a number $M =M (c_V ) <\infty$, such that;
$$
\max_{\phi}{\cal E}^B_{j,\phi}\int_0^{\tau_B} 
{\Bbb I}_{\{X(s)=j\}}\dr s~\leq~M^k,
$$
whenever $\text{d} (j,B) =k$. On the other hand, by Theorem~\ref{Main},
$$
\sum_{B:\,\text{d}(j,B)=k}\rho (B)~\leq~\e^{-\nu k^2} ,
$$
as soon as $k>R$. Therefore, the sum in \eqref{selfterm} converges, and the
proof of Lemma~\ref{Lemma1} is, thereby, concluded

\qed
\vskip 0.2in
\noindent
{\it Proof of Theorem~\ref{Main}:}\ 
Let us start by introducing some additional notation: Given a finite set
$B\subset\Ztwo$ with the decomposition \eqref{Bdecomp} into the disjoint union of
connected components  $B_1,...,B_n$ we say that another set $A$ is a dry
neighbour of $B$; $A\in{\cal D}_B$, if
$$
A\cap B=\emptyset\qquad\text{but}\qquad D\cup\partial 
B_l\neq\emptyset
;\ l=1,...,n.
$$
\begin{prop}
\label{Prop3}
There exists a constant $c_4 =c_4 (J)$, such that for every finite 
$B\subset\Ztwo$,
\begin{equation}
\label{dryneighbour}
\sum_{A\in{\cal D}_B}\rho (A)~\leq~\e^{-c_4 |B|} .
\end{equation}
\end{prop}

The proof of Proposition~\ref{Prop3} relies on the following two  basic
estimates which have been proven in \cite{DV}:
\begin{enumerate}
\item There exists a number $M=M(J)$ and a constant $c_5 =c_5 (J)$, such that,
\begin{equation}
\label{Yvan1}
\inf_{A\in{\cal D}_B}\sum_{C\subset B}\e^{J|C|}\frac{{\bf Z}(A\cup 
C)}
{{\bf Z}(A)}~\geq~\e^{c_5 |B|} ,
\end{equation}
whenever $B$ is connected and $\text{diam}(B)\geq M$.
\item Let $A\neq\emptyset$ and $i\in\Ztwo\setminus A$. Then,
\begin{equation}
\label{Yvan2}
\frac{{\bf Z}(A\cup\{i\})}{{\bf Z}(A)}~\geq~\frac{c_6 (J)}{\sqrt{\dr
    (i,A)}} .
\end{equation}
\end{enumerate}

The above estimates are linked to the claim of  Proposition~\ref{Prop3} in the
following way:
$$
\sum_{A\in{\cal D}_B}\rho (A)
~\leq~\left(\inf_{A\in{\cal D}_B}\sum_{C_1\subset
  B_1}...\sum_{C_n\subset
B_n}\frac{{\bf Z}(A\cup_1^nC_l)}{{\bf
  Z}(A)}\e^{J\sum_1^n|C_l|}\right)^{-1}.
$$
If, for some $m\in [1,...,n-1]$, we regroup $B$ as 
$$
B~=~B^+\cup B^-~\df~\left\{ B_1
,...,B_m\right\}\bigcup\left\{B_{m+1},...,
B_n\right\} ,
$$
then,  since $A\cup_1^m C_l$ always belongs to ${\cal D}_{\cup_{m+1}^n B_l}$,
we obtain the following decoupling estimate:
\begin{equation}
\label{decouple}
\begin{split}
\inf_{A\in{\cal D}_B}&\sum_{C_1\subset
  B_1}...\sum_{C_n\subset
B_n}\frac{{\bf Z}(A\cup_1^nC_l)}{{\bf
  Z}(A)}\e^{J\sum_1^n|C_l|}\\
&\geq~\inf_{A\in{\cal D}_B^+}\sum_{C_1\subset
  B_1}...\sum_{C_m\subset
B_m}\frac{{\bf Z}(A\cup_1^m C_l)}{{\bf
  Z}(A)}\e^{J\sum_1^m|C_l|}\\
&\qquad\times\inf_{A\in{\cal D}_B^-}\sum_{C_{m+1}\subset
  B_{m+1}}...\sum_{C_n\subset
B_n}\frac{{\bf Z}(A\cup_{m+1}^nC_l)}{{\bf
  Z}(A)}\e^{J\sum_{m+1}^n|C_l|}  .
\end{split}
\end{equation}

In particular,  the claim \eqref{dryneighbour} directly follows from the 
estimate \eqref{Yvan1} whenever $\text{diam}(B_l ) >M$ for each $l=1,...,n$. In
fact, in view of \eqref{Yvan1} and \eqref{decouple},   it remains to study only
the case when all connected components of  $B$ are small;  $\text{diam}(B_l )
<M;\ l=1,...,n$.

In the latter situation, however,  we can use  \eqref{Yvan2} and  estimate;
$$
\frac{{\bf Z}(A\cup C_l)}{{\bf Z}(A)}~\geq~
\left(\frac{c_6}{\sqrt{2M+1}}
\right)^{|C_l |},
$$
for every $l$; $A\in{\cal D}_{B_l}$ and $C_l\subset B_l$. Therefore,
$$
\inf_{A\in{\cal D}_B}\sum_{C_1\subset
  B_1}...\sum_{C_n\subset
B_n}\frac{{\bf Z}(A\cup_1^nC_l)}{{\bf
  Z}(A)}\e^{J\sum_{l=1}^n|C_l|}~
\geq~\prod_1^n\left(1+\frac{c_6}{\sqrt{2M+1}}\right)^{|B_l |}
,
$$
and \eqref{dryneighbour} follows.

\qed
\vskip 0.1in

\begin{rem}
One could hope to deduce from Proposition~\ref{Prop3} the claim of 
Theorem~\ref{Main} even without the additional assumption \eqref{diam}. We were
not able to do so, and, moreover, even not sure that the corresponding
statement would be true --- the entropy cancelation forced by the  condition
\eqref{diam} could well be essential for the validity of  the claim. We would
like to stress, however, that within the framework of the renormalization
approach we try to develop there is absolutely no point to relax \eqref{diam}.
\end{rem}

\noindent The rest of the proof is an adaptation of the  ideas of \cite{DV} to
the case of multiply connected sets:
\vskip 0.1in
First of all, for any finite $D\subset\Ztwo$ let us denote its  $k$-enlargement
$D^{(k)}$ as
$$
D^{(k)}~\df~\left\{ i\in\Ztwo :\ \dr (i,D )\leq k\right\} .
$$
Assume now that $B=\bigvee_1^n B_l$ is as in  the assumptions of 
Theorem~\ref{Main}, that is the  diameter of each connected component $B_l$ of
$B$ is bounded below,   $\text{diam} (B_l )\geq R;\ i=1,...,n$. 

We have to show that the bound \eqref{mainbound} holds uniformly in such $B$-s
as soon as $R$ is chosen large enough.

Let us say that a tuple $\underline{k} =(k_1 ,...,k_n )$  of $n$ natural
numbers is $B$-admissible if:
\begin{itemize}
\item $k_1\in{\Bbb N}$ (no restriction).
\item either $k_2 =0$, or the sets $B^{(k_1 )}_1$ and $B_2^{(k_2)}$ 
are
disjoint.
\item either $k_3 =0$, or the set $B_3^{(k_3)}$ is disjoint from
$$
B_1^{(k_1)}\cup B_2^{(k_2)} .
$$
\item ...................................
\item either $k_n=0$, or the set $B_n^{(k_n)}$ is disjoint from
$$
\bigcup_1^{n-1}B_l^{(k_l)} .
$$
\end{itemize}
For any $B$-admissible tuple $\underline{k}$ we set
$$
B^{(\underline{k})}~\df~\bigcup_1^{n}B_l^{(k_l)} .
$$
This construction enjoys the following two properties:
\begin{enumerate}
\item For any $A\cap B =\emptyset$ there is the unique $B$-admissible
tuple $\underline{k}$, such that,
$$
A\in{\cal D}_{B^{(\underline{k})}} .
$$
Indeed, this tuple $\underline{k} =(k_1 ,...,k_n)$ can be constructed in the
following way:
\begin{equation*}
\begin{split}
&~k_1~=~\max \big\{ k:\ B_1^{(k)}\cap A =\emptyset\}\\
&~k_2~=~\max \big\{ k>0:\ B_2^{(k)}\cap (A\cup B_1^{(k_1)})
=\emptyset\} \\ 
&~\cdot\\
&~\cdot\\
&~\cdot\\
&~k_n~=~\max \big\{ k>0:\ B_n^{(k)}\cap (A\cup_1^{n-1} B_l^{(k_l)}) 
=\emptyset\} \\ 
\end{split}
\end{equation*}
with the convention that the maximum over an empty set equals zero.
\item For any $B$-admissible tuple $\underline{k}=(k_1,...,k_n)$;
$$
\left| B^{(\underline{k} )}\right|~\geq~|B| +\sum_1^n k_l .
$$
This follows directly from the definition of the $B$-admissibility.
\end{enumerate}

Using Proposition~\ref{Prop3} we, thereby, obtain:
\begin{equation*}
\begin{split}
\sum_{A\cap B =\emptyset}\rho
(A)~&=~\sum_{B-\text{admissible}\,\underline{k}}
\sum_{A\in{\cal D}_B^{(\underline{k})}}\rho (A)\\
&\leq~\sum_{B-\text{admissible}\,\underline{k}}\e^{-c_4 (|B|+\sum
  k_l)}\\
&\leq~\e^{-c_4 |B|}\left( 1- \e^{-c_4}\right)^{-n} .
\end{split}
\end{equation*}
By the assumption \eqref{diam}, $n\leq |B|/R$. Thus it remains to choose $R
=R(J)$ so large that,
$$
\nu (J)~\df~c_4 (J) +\frac{\log (1-\e^{-c_4 (J)})}{R}~>~0,
$$
and \eqref{mainbound} follows.
\qed

\end{document}